\documentclass{amsart}
\usepackage[mathscr]{eucal}
\usepackage{latexsym}
\usepackage{amsmath}
\usepackage{amssymb}
\usepackage{amsthm}
\usepackage{graphicx}

\newtheorem{theorem}{Theorem}[section]
\newtheorem{lemma}[theorem]{Lemma}
\newtheorem{corollary}[theorem]{Corollary}
\theoremstyle{definition}
\newtheorem{example}[theorem]{Example}
\theoremstyle{remark}
\newtheorem{remark}[theorem]{Remark}
\numberwithin{equation}{section}

\newcommand{\im}{{\rm Im}\hspace{0.02in}}
\newcommand{\isom}{{\rm Isom}\hspace{0.02in}}

\title{Global Homeomorphisms and Covering Projections on Metric spaces}
\author{Olivia Gut{\'u}}
\address{Centro de Investigaci\'{o}n en Matem\'{a}ticas, Universidad Aut\'{o}noma de Hidalgo, 42184 Pachuca, Mexico}
\email{olivia@uaeh.edu.mx}
\thanks{The first author was supported in part  by {\sc Promep} (M{\'e}xico) Grant
103.5/03/2568.}

\author{Jes{\'u}s A. Jaramillo}
\address{Departamento de An{\'a}lisis Matem{\'a}tico, Universidad Complutense de Madrid, 28040 Madrid, Spain}
\email{jaramil@mat.ucm.es}
\thanks{Both authors were  supported in part by {\sc D.G.E.S.} (Spain)
Grant BFM2003-06420.}

\subjclass{58C15, 58B20, 46T05}

\begin{document}

\maketitle

\begin{abstract}
For a large class of metric spaces with nice local structure, which
includes Banach-Finsler manifolds and geodesic spaces of curvature
bounded above, we give sufficient conditions for a local
homeomorphism to be a covering projection. We first obtain a general
condition in terms of a path continuation property. As a
consequence, we deduce several conditions in terms of path-liftings
involving a generalized derivative,  and in particular we obtain an
extension of Hadamard global inversion theorem in this context. Next
we prove that, in the case of quasi-isometric mappings, some of
these sufficient conditions are also necessary. Finally, we give
some applications to the existence of global implicit functions.
\end{abstract}

\section{Introduction}

Let $f:{\mathbb R^n} \to {\mathbb R^n}$ be a $C^1$ map with
everywhere nonvanishing Jacobian. A natural question is to ask under
which conditions we can assure that $f$ is a global diffeomorphism
(or, equivalently, a global homeomorphism). This problem was first
considered by Hadamard \cite{hadamard}, who obtained a sufficient
condition in terms of the growth of $\Vert [df(x)]^{-1}\Vert$, by
means of his celebrated  {\it integral condition}. Namely, $f$ is a
global diffeomorphism provided
$$
\int _0^\infty \inf_{\vert x \vert = t}\|[df(x)]^{-1}\|^{-1} ~
dt=\infty.
$$
This result was  extended to the infinite-dimensional setting by P.
L{\'e}vy \cite{levy}, who considered the case of a $C^1$ mapping $f$
between Hilbert spaces.

Some years later, Cacciopoli \cite{cacciopoli} and Banach and Mazur
\cite{banach} obtained a purely topological condition of global
inversion for a local homeomorphism $f$ between Banach spaces.
Namely, they proved that $f$ is a global homeomorphism if, and only
if, $f$ is a proper map. Properness was also used by Palais
\cite[Section 4]{palais0}, in the context of locally compact spaces
and (finite-dimensional) manifolds. The properness condition was
relaxed to closedness by Browder \cite{browder} in a more general
setting.

Later on, Plastock  extended in \cite{plastock} the
Hadamard-L{\'e}vy theorem to the case of $C^1$ mappings between
Banach spaces. In fact, Plastock obtained a more general result as
follows. He introduced a limiting property for lines, called  {\it
condition (L)}, which is analogous to the continuation property used
by Rheinboldt in \cite{rheinboldt}. Then, for a local homeomorphism
between Banach spaces, Plastock proved that $f$ satisfies condition
(L) if, and only if, $f$ is a covering projection (or, equivalently,
a global homeomorphism). Finally, he showed that properness,
closedness or the Hadamard integral condition all imply condition
(L). Since then, this condition (L) has proved to be quite useful in
global inversion theorems, and it has been widely used, as can be
seen for instance in \cite{radulescu}, \cite{pourciau1},
\cite{ioffe} or \cite{partha}.

The question of global invertibility has been also studied from the
point of view of nonsmooth analysis. In this sense, the problem of
giving analytic conditions of global inversion for a nonsmooth
mapping $f$ between Banach spaces was first considered by F. John in
\cite{john}. In particular, John obtained an extension of the
Hadamard integral condition in this setting, in terms of the lower
scalar Dini derivative of $f$. For the proof, he used the
prolongation of local inverses of $f$ along lines. Further results
in this direction were given by Ioffe \cite{ioffe} in terms of the
so-called {\it surjection constant} of the mapping $f$, making use
of the aforementioned condition (L). In the finite-dimensional case,
analogous results were obtained by Pourciau \cite{pourciau1},
\cite{pourciau2} by means of Clarke generalized jacobian of $f$. The
surjection constant of a mapping $f$ was also considered by  Katriel
\cite{katriel} in order to obtain global inversion theorems in
certain metric spaces. The methods of Katriel came from critical
point theory, and in particular are based on a suitable Mountain
Pass theorem.

Global inversion problems have been less studied  in the context of
smooth mappings between Riemannian or Finsler manifolds. In the
Riemannian case, a classical result due to Ambrose \cite{ambrose}
gives conditions for a  local isometry to be a global
diffeomorphism. This was extended by Wolf and Griffiths \cite{WG},
who obtained more general conditions under which a local
diffeomorphism is a covering projection. On the other hand, we have
a quite general result due to Rabier \cite{rabier} for Finsler
manifolds, which establishes the global inversion of $f$ using a
growth condition for $\Vert [df(x)]^{-1}\Vert$ that does not require
integrals.

Our aim in this paper is twofold. On one hand, we provide an
extension of the above mentioned results to the framework of metric
spaces. On the other hand, we present them in a unified and
systematic way, where the ideas and methods of line-lifting play a
central role, which leads to a clarification and simplification of
proofs. To this end, we introduce in Section 2 a fairly general
class of path-connected metric spaces with nice local structure,
namely the class of metric spaces which are $\mathcal P$-connected
and locally $\mathcal P$-contractible spaces, which we will define
below. These include Banach spaces and Banach manifolds,  as well as
many other ``singular" spaces, as for example geodesic metric spaces
of curvature bounded above. If $X$ and $Y$ belong to this class of
spaces, our main goal is to find conditions for a local
homeomorphism $f:X \to Y$ to be a covering projection. We first
obtain in Theorem \ref{teorema1} a general condition, in terms of a
continuation property. This Theorem is the key of our presentation,
since every further result will be derived from it. Next, in order
to give analytical conditions in this nonsmooth setting, we consider
in Section 3 the upper and lower scalar Dini derivatives of the
mapping $f$ at $x$, denoted respectively by $D_x^+ f$ and $D_x^- f$
(in the smooth case, these quantities reduce to $\Vert df(x)\Vert$
and $\Vert [df(x)]^{-1}\Vert^{-1}$). Then we obtain in Theorem
\ref{teorema2} and Theorem \ref{teorema3} two mean value
inequalities in this context, which are going to be quite useful in
the sequel. In Section 4 we introduce a bounded path-lifting
condition in terms of $D_x^- f$, and in Theorem \ref{teorema4} we
see that it is a sufficient condition for $f$ to be a covering
projection. We derive some consequences, and in particular we obtain
a version of Hadamard integral condition (see Theorem
\ref{teorema5}) in our setting. In Section 5 we provide a more
complete result under some extra regularity conditions on $f$. More
precisely, for a quasi-isometric mapping $f$ we give  in Theorem
\ref{teorema6} several conditions which are necessary and sufficient
for $f$ to be a covering projection or a global homeomorphism.
Finally, Section 6 is devoted to applying our results in order to
obtain global implicit function theorems in certain metric spaces.

\section{Continuation property on metric spaces}

Our purpose in this Section is to give a general condition for a
local homeomorphism between metric spaces to be a covering
projection. This will be achieved by means of a continuation
property, much in the spirit of Rheinboldt \cite{rheinboldt} and
Plastock \cite{plastock}. Our result will apply to a wide class of
path-connected metric spaces, which we introduce now.

\

Let  $Y$ be a  metric space, and let $\mathcal P$ be a family of
continuous paths in $Y$. We say that $Y$ is {\it $\mathcal
P$-connected} if the following conditions hold:

\begin{enumerate}
\item If the path $p:[a,b]\rightarrow Y$ belongs to $\mathcal P$, then the
reverse path $\bar p$, defined by $\bar p(t)=p(a-t+b)$, also belongs
to $\mathcal P$.
\item Every two points
in $Y$ can be joined by a path in $\mathcal P$.
\end{enumerate}
We say that $Y$ is {\it locally $\mathcal P$-contractible} if
every point $y_0\in Y$ has an open neighborhood $U$ which is {\it
$\mathcal P$-contractible}, in the sense that there exists a
homotopy $H:U\times [0,1]\rightarrow U$ satisfying:
\begin{enumerate}
\item [(3a)] $H(y_0,t)=y_0$,  for all $t \in [0,1]$.
\item [(3b)] $H(y,0)=y_0$ and $H(y,1)=y$, for all $y\in U$.
\item [(3c)]  For every $y\in U$, the path $p_y(t):=H(y,t)$ belongs to $\mathcal P$.
\end{enumerate}
Next we give some general  examples of spaces satisfying these
conditions.
\begin{example}
It is clear that every normed vector space $V$ is $\mathcal
L$-connected and locally $\mathcal L$-contractible, where
$\mathcal L$ is the family of all lines in $V$, that is, paths of
the form $p(t)=(1-t)y_0+ty_1$, with $y_0, y_1\in V$. The same is
true for any convex subset of $V$.
\end{example}

\begin{example}\label{ejemplo2}
Let $M$ be a connected paracompact Banach manifold of class $C^k$,
for $0\leq k\leq \infty$, and let $\mathcal P^k$ denote the family
 of all $C^k$-paths on $M$. Since $M$ is paracompact, then it is
metrizable (see \cite{palais}). It is easy to see that, with any
equivalent metric, $M$ is $\mathcal P^k$-connected and locally
$\mathcal P^k$-contractible. Indeed, the reverse path of a
$C^k$-path in $M$ is itself $C^k$. By connectedness, every two
points in $M$ can be joined by a $C^k$-path. Furthermore, for
every $y_0\in M,$ a  $\mathcal P^k$-contractible neighborhood can
be obtained by considering $U=\phi^{-1}_\alpha(B_r(u_0))$, where
$(U_\alpha,\phi_\alpha)$ is a $C^k$ chart of $M$ about $y_0$ and
$B_r(u_0)$ denotes the ball of radius $r$ centered at
$u_0=\phi_\alpha(y_0)$, where $r>0$ is small enough such that
$B_r(u_0)\subset\phi_\alpha(U_\alpha)$. Then define $H:U\times
[0,1]\rightarrow U$ by
$H(y,t)=\phi_\alpha^{-1}((1-t)\phi_\alpha(y_0)+t\phi_\alpha(y))$.
\end{example}

\begin{example}
Let $M$ be a connected  $n$-dimensional Lipschitz manifold (see e.
g. \cite{LV} or \cite{CP}) and suppose that $M$ is endowed with a
metric which is locally Lipschitz equivalent to the Euclidean one.
If $\mathcal P_L$ denotes the family of all Lipschitz paths in
$M$, it can be seen as before that $M$ is $\mathcal P_L$-connected
and locally $\mathcal P_L$-contractible.
\end{example}

Along the paper, we will focus on the family  $\mathcal R $ of all
{\it rectifiable} paths on a given metric space, and we will
consider accordingly the class of {\it locally  $\mathcal R
$-contractible} spaces. Recall that, for a metric space $Y$, the
{\it length} of a path $p:[a,b]\rightarrow Y$ is defined by:
$$\ell(p)=\sup_{a=t_0\leq t_1 \leq \ldots \leq t_n=b}\sum_{i=0}^{n-1}d(q(t_i),q(t_{i+1})),$$
where the supremum is taken over all partitions  $a=t_0\leq t_1
\leq \ldots \leq t_n=b$ (no bound on $n$). The path $p$ is  said
to be {\it rectifiable} when $\ell(p)<\infty$.  It is not
difficult to see that if $Y$ is path-connected and locally
$\mathcal R $-contractible then $Y$ is also $\mathcal R
$-connected.

\begin{example} Now we describe some classes of metric spaces which  are
 locally  $\mathcal R$-contractible.
\begin{enumerate}
\item {\it Normed spaces.} Of course every convex subset of a
normed vector space is $\mathcal R$-connected and  locally
$\mathcal R$-contractible.

\item  {\it Finsler manifolds.} Let $E$ be a real Banach space and let $M$
be an $E$-manifold of class $C^1$, with atlas
$\{(U_\alpha,\phi_\alpha):\alpha\in\Lambda\}$ and tangent bundle
$TM$. For every $\alpha$ we can identify $TU_\alpha$ with
$U_\alpha\times E$.  Recall that a functional
$\|\cdot\|:TM\rightarrow \mathbb{R}$ is a {\it Finsler structure}
for $TM$ (according to \cite{palais2}) if, for every
$\alpha\in\Lambda$, the following conditions hold:
\begin{enumerate}
\item for every $x\in U_\alpha$ the map $\|\cdot\|_x$,  defined by $v\mapsto \|(x,v)\|$, is
an admissible norm for $E$.
\item for every $x_0\in U_\alpha$ and $s>1$, there exists a
neighborhood $U_{\alpha(s)}$ of $x_0$ such that:
$$s^{-1}\|v\|_{x_0}\leq\|v\|_x\leq s\|v\|_{x_0}$$
for all $x\in U_{\alpha(s)}$ and $v\in E$.
\end{enumerate}
If there exists such a functional for $TM$, we say that $M$ is a
{\it Finsler manifold}. In this case, the  {\it Finsler length} of a
$C^1$ path $p:[a,b]\rightarrow M$ is defined as:
$$\ell_F(p)=\int_a^b\| p'(t)\|~dt.$$
On other hand, if $M$ is  connected,   for every $x,y\in M$, there
exists a $C^1$ path  joining $x$ to $y$, and we can define the {\it
Finsler distance} by:
$$d_F(x,y)=\inf\{\ell_F(p): p\text{ is a  $C^1$ path from } x \text{
to } y \}.$$  We always consider that $M$ is endowed with this
metric, which is compatible with the topology of of $M$ (see
\cite{palais2}). Therefore, every $C^1$-path in $M$ is
rectifiable. As a consequence, we obtain as in Example
\ref{ejemplo2} that every connected $C^1$ Finsler manifold is
$\mathcal R$-connected and locally $\mathcal R$-contractible. In
particular this includes connected Riemannian manifolds, both in
the finite-dimensional and infinite-dimensional cases (see e. g.
\cite{lang}).

\item {\it Lipschitz manifolds.} Every Lipschitz path in a metric
space is rectifiable. Therefore if $M$ is a connected
$n$-dimensional Lipschitz manifold  endowed with a metric which is
locally Lipschitz equivalent to the Euclidean one, then $M$ is
$\mathcal R$-connected and locally  $\mathcal R$-contractible.

\item {\it Geodesic spaces}. Recall that a path $g:[0,1]\rightarrow Y$ in
a metric space $Y$ is said to be a (constant speed) {\it geodesic}
if there exists $L>0$ such that  $d(g(t),g(t'))=L \cdot|t-t'|$, for
all $t,t'\in[0,1]$. Note that in this case $g$ is rectifiable and
$\ell(g)=L$. We say that $Y$ is a {\it geodesic space} if  every two
points in $Y$ can be joined by a geodesic.

Now suppose that $Y$ is a geodesic space such that every point
$y_0\in Y$ has an open neighborhood $U$ verifying: (a) for each
$y\in U$ there exists a unique geodesic $g_y$ in $U$ from $y_0$ to
$y$, and (b) if $y_n\rightarrow y$ in $U$, then
$g_{y_n}(t)\rightarrow g_y(t)$ uniformly on $[0,1]$. Then we can
apply Lemma \ref{lema1} below to the  mapping
$H:U\times[0,1]\rightarrow U$ given by $H(y,t)=g_y(t)$ and  we
obtain that $Y$ is $\mathcal R$-connected and  locally  $\mathcal
R$-contractible.

In the book by Bridson and Haefliger \cite{bridson} we can find
several classes of geodesic spaces satisfying the above
requirements. For example, this is the case of proper geodesic
spaces which are locally uniquely geodesic (see
\cite[I.3.13]{bridson}). It is also the case of geodesic spaces of
curvature $\leq \kappa$ (see \cite[II.1.4]{bridson}). Spaces of this
kind  include a large class of polyhedral complexes (see
\cite[II.5.5]{bridson}; see also \cite[I.7.57]{bridson} and the
comments before \cite[I.7.57]{bridson}).
\end{enumerate}
\end{example}

\begin{lemma}\label{lema1}
Let $U$ be an open set  in a metric space $Y$, consider a mapping
$H:U\times[0,1]\rightarrow U$, and set  $p_y(t):=H(y,t)$. The
following statements are equivalent:

\begin{enumerate}
  \item The map $H$  is continuous
  on $U\times [0,1]$.
  \item For every $y\in U$ and $\varepsilon>0$, there exists $\delta>0$ such that
  if $z\in B_\delta(y)\subset U$ then $$d(p_z(t),p_y(t))<\varepsilon, \hspace{0.1in}\forall t\in[0,1].$$
  \item For all $y\in U$, if $y_n\rightarrow
y$ in $U$, then $p_{y_n}(t)\rightarrow p_y(t)$ uniformly on $[0,1]$.
In this case, we say that the paths $\{p_y:y\in U\}$ vary
continuously with their endpoints
\end{enumerate}
\end{lemma}

\begin{proof}
$(1 \Rightarrow 2)$ Suppose that $H$ is continuous on $U\times
[0,1]$. Then, for every $y\in U$, $H$ is uniformly continuous over
the compact $K=\{y\}\times [0,1]$. Given $\varepsilon
>0$, there exists $\rho>0$ such that, if $s,t\in[0,1]$, with $|s-t|<\rho$,
then $d(p_y(s),p_y(t))<\frac{\varepsilon}{2}.$ On the other hand,
for all $t\in[0,1]$ we can take $0<\delta_t<\rho$ such that, if
$z\in B_{\delta_t}(y)$ and $s\in (t-\delta_t,t+\delta_t)$, then
    $d(p_z(s),p_y(t))<\frac{\varepsilon}{2}.$
By the compactness of  $[0,1]$, there exist $t_1, ... , t_m\in
[0,1]$, such that
    $[0,1]\subset\bigcup_{j=1}^m(t_j-\delta_{t_j}, t_j+\delta_{t_j}).$
Consider $\delta=\min\{\delta_{t_1}, ... ,\delta_{t_m}\}>0$. Let
$z\in B_\delta(y)$. Then, for  $t\in [0,1]$ there exists
$j\in\{1,...,m\}$ such that if $t\in
(t_j-\delta_{t_j},t_j+\delta_{t_j})$, then
$d(p_z(t),p_y(t_j))<\frac{\varepsilon}{2}$ and
$d(p_y(t_j),p_y(t))<\frac{\varepsilon}{2}$, because
$|t-t_j|<\delta_j<\rho$. Therefore, $d(p_z(t),p_y(t))<\varepsilon.$

\noindent $(2 \Rightarrow 3)$  It is straightforward.

\noindent $(3 \Rightarrow 1)$ Since $p_y$ is continuous, if
$t_n\rightarrow t$ we have that $d(p_y(t_n),p_y(t))\rightarrow 0$.
If $y_n\rightarrow y$, by hypothesis,
$d(p_{y_n}(t_n),p_y(t_n))\rightarrow 0$. Then,
$d(p_{y_n}(t_n),p_y(t))\rightarrow 0$.
\end{proof}

Let $X$ and $Y$ be metric spaces, and let $p:[0,1]\rightarrow Y$ be
a path in $Y$. We will say that a continuous map $f:X\rightarrow Y$
has the {\it continuation property for $p$} if, for every $b\in
(0,1]$ and every continuous path $q:[0,b)\rightarrow X$ such that
$f\circ q=p$ over $[0,b)$, there exists a sequence $\{t_{n}\}$ in
$[0,b)$ convergent to $b$ and such that $\{q(t_{n})\}$ converges in
$X$.

\begin{theorem}\label{teorema1}
Let $f:X\rightarrow Y$ be a local homeomorphism between
 metric spaces and suppose that $Y$ is  $\mathcal P $-connected and
 locally $\mathcal P $-contractible for some family $\mathcal P $ of paths.
Then  $f$ is a covering projection if and only if $f$ has the
continuation property for every path in $\mathcal P$.
\end{theorem}

\begin{proof}
The sufficiency follows from  general properties of  covering
projections. If $f$ is a covering projection, then $f$ lifts paths
and has the unique-path-lifting property (see for example
\cite[Section 2.2]{spanier}). That is, for every path
$p:[0,1]\rightarrow Y$ in $\mathcal P$, there exists a  unique  path
$q_0:[0,1]\rightarrow X$ such that $f\circ q_0=p$ over $[0,1]$. Thus
every partial lifting $q:[0,b)\rightarrow X$ of $p$ with $0<b\leq 1$
must coincide with $q_0$ over $[0,b)$, and therefore $f$ has the
continuation property for $p$.

For the necessity, suppose that $f$ has the continuation property
for each path in $\mathcal P$. We first prove that $f$ lifts every
path $p:[0,1]\rightarrow Y$ in  $\mathcal P$. Indeed, let
$y=p(0)\in f(X)$ and consider $x\in f^{-1}(y)$. Since $f$ is a
local homeomorphism, there exist $\varepsilon >0$ and a path
$q:[0,\varepsilon)\rightarrow X$ beginning at $x$ such that
$f\circ q=p$ over $[0,\varepsilon)$. Consider  the largest number
$b$ in $[0,1]$ for which $q$ can be extended to a continuous path
such that $f\circ q=p$ over $[0,b)$, and suppose that $b<1$. There
exists a sequence $\{t_{n}\}$ in $[0,b)$ convergent to $b$ and
such that $\{q(t_{n})\}$ converges to some point $z\in X$. By
continuity $f(z)=p(b)$. Now let $V$ be a neighborhood of $z$ such
that $f|_V$ is a homeomorphism. Then there exists some $n_0\in
\mathbb N$ such that $q(t_{n})\in V$ for $n\geq n_0$. Also, there
exists $\delta>0$ and a path $q_{1}:(b-\delta,b+\delta)\rightarrow
X$ with $q_{1}(t_{m})=q(t_{m})$, where $m\geq n_0$ and $b -\delta
<t_{m}<b$, and such that $f\circ q_{1}=p$ over $(b -\delta,b
+\delta)$. Since local homeomorphisms have the unique-lifting
property, $q$ can be extended to a continuous path (call it again
$q$) over $[0,b +\delta)$, beginning at $x$ and such that $f\circ
q=p$ over $[0,b +\delta)$. This contradicts the maximality of
$[0,b)$. Therefore $b=1$ and the same reasoning shows that $q$ can
be extended to $[0,1]$.

Now let $y\in Y,$ $y'\in f(X)$ and $p:[0,1]\rightarrow Y$ be a path
in $\mathcal P$ joining $y'$ to $y$. Then there exists $q$ such that
$f\circ q=p$ over $[0,1]$; in particular $f(q(1))=p(1)=y$.
Therefore, $f$ is onto. Next we are going to show that $f$ is a
covering projection.

Let $y_0\in Y,$ $x\in f^{-1}(y_0)$ and let $U$ be a $\mathcal P
$-contractible neighborhood of $y_0$. Let $q_y$ be the unique
lifting of $p_y$, with $q_y(0)=x$ and such that $f\circ q_y=p_y$.
Now consider:
 $$O_x=\{q_y(1):y\in U\}.$$

We will prove that $f|_{O_x}:O_x\rightarrow U$ is a homeomorphism
and  that $f^{-1}(U)$ is the disjoint union of open sets $O_x$,
$x\in f^{-1}(y_0)$.

It is easy to see that  $f|_{O_x}:O_x\rightarrow U$ is bijective.
Since $f$ is a local homeomorphism, in order to see that $f|_{O_x}$
is a homeomorphism we just have to prove that $O_x$ is open in $X$.
Let $x'\in O_x$ and $y\in U$ such that $x'=q_y(1)$. For every
$u\in\im q_y$ there exist an open neighborhood  $U^u$  and an open
ball $B_{r_u}(f(u))$, such that $f|_{U^u}:U^u\rightarrow
B_{r_u}(f(u))$ is a homeomorphism. Let $V^u\subset U^u$ be an open
set such that $f(V^u)=B_{\frac{r_u}{2}}(f(u))$. By compactness,
there exist $u_1,...,u_m\in \im q_y$ such that $\im q_y\subset
V^{u_1}\cup ... \cup V^{u_m}.$ For $k=1,2,...,m$, let us denote
$V_k=V^{u_k}$, $U_k=U^{u_k}$, $y_k=f(u_k)$ and $r_k=r_{u_k}$. Then,
 $$\im p_y\subset \bigcup_{k=1}^m B_{\frac{r_k}{2}}(y_k)=\bigcup_{k=1}^m B_k.$$
Let $s_k:B_{r_k}(y_k)\rightarrow U_k$ be the inverse of
$f|_{U_k}$. Let $\rho>0$ be the Lebesgue's number of $[0,1]$ for
the finite covering $\left\{q_y^{-1}(V_k\cap\im
q_y)\right\}_{k=1}^m$ and let $0=t_0<t_1<\ldots<t_n=1$ be a
partition of $[0,1]$ such that, for  every $j=1,\ldots,n$,  the
diameter of $[t_{j-1},t_j]$ is less than $\rho$. Then for each
$j=1,\ldots,n$ there exists $k_j\in\{1,\ldots,m\}$ such that
$q_y[t_{j-1},t_j]\subset V_{k_j}$. For each $j=1,\ldots, n$, let
$\tilde u_j=q_y(t_j)\in V_{k_j}\cap V_{k_{j+1}}=\tilde V_j$. Since
$\tilde V_j$ is open, then $f(\tilde V_j)$ is  open in $Y$ and
contains $\tilde y_j=p_y(t_j)$. Also, $f|_{\tilde V_j}:\tilde
V_j\rightarrow f(\tilde V_j)$ is a homeomorphism and
$s_{k_j}\equiv s_{k_{j+1}}$ over $f(\tilde V_j)$. Let $\delta_j>0$
be such that $B_{\delta_j}(\tilde y_j)\subset f(\tilde V_j)$; in
particular $s_{k_j}\equiv s_{k_{j+1}}$ over
 $B_{\delta_j}(\tilde y_j)$. Let $\varepsilon>0$ be such that
    $$0<\varepsilon<{\rm dist}\hspace{0.02in}(\im p_y ; Y\setminus (B_1\cup\ldots \cup B_m))$$
and $0<\varepsilon<\min\{r_1,\ldots,r_m;\delta_1,...,\delta_n\}$.
Therefore, there exists $\delta>0$ such that, if $z\in B_\delta(y)$
and $B_\delta(y)\subset B_{r_{k_n}}(y_{k_n})$, then
$d(p_z(t),p_y(t))<\frac{\varepsilon}{2}$, $\forall t\in[0,1].$ For
all $j=1,\ldots,n$; if $t\in[t_{j-1},t_j]$ then
$d(p_y(t),y_{k_j})<\frac{r_{k_j}}{2}$, therefore:
    $$d(p_z(t),y_{k_j})\leq d(p_z(t),p_y(t))
    +d(p_y(t),y_{k_j})<\frac{\varepsilon}{2}+
    \frac{r_{k_j}}{2}<r_{k_j}.$$
In other words, $p_z[t_{j-1},t_j]\subset B_{r_{k_j}}(y_{k_j})$
where the local inverse $s_{k_j}$ is defined. Furthermore, for all
$j=1,\ldots,n$:
    $$d(p_z(t_j),p_y(t_j))=d(p_z(t_j),\tilde y_j)<\varepsilon<\delta_j$$
so, $p_z(t_j)\in B_{\delta_j}(\tilde y_j)$ where we know that
$s_{k_j}\equiv s_{k_{j+1}}$. Therefore,
    $s_{k_{j}}(p_z(t_j))=s_{k_{j+1}}(p_z(t_j))$,  for all $j=1,\ldots,n.$
Then, the lifting $q_z$ of $p_z$ with $q_z(0)=x$ satisfies
$q_z(t)=s_{k_j}\circ p_z(t)$,  for all $t\in [t_{j-1},t_j].$
Therefore, the set $s_{k_n}(B_\delta(y))$ is an open set
containing $x'=q_y(1)$ and contained in $O_x$.

We will prove that $f^{-1}(U)=\bigcup_{x\in f^{-1}(y_0)}O_x$. Let
$x\in f^{-1}(y_0)$ and $x'\in O_x$. Then $x'=q_y(1)$ for some
$y\in U$. Therefore $f(x')=p_y(1)=y$, so $x'\in f^{-1}(U)$. On the
other hand, let $x'\in f^{-1}(U)$; there exists $y\in U$ such that
$f(x')=y$. Let $\bar p(t)=p_y(1-t)$. Since $f$ lifts paths in
$\mathcal P$, there exists $\bar q$ such that $f\circ\bar q=\bar
p$ with $\bar q(0)=x'$. Setting $q_y(t)=\bar q(1-t)$ we get
$f\circ q_y=p_y$, with $q_y(1)=x'$. Therefore $x'\in O_{q_y(0)}$.

Let $x'\in O_{x_1}\cap O_{x_2}$. For $i=1,2$, there exists $y_i\in
U$, the paths $p_{y_i}$ and their liftings $q_{y_i}$ such that
$q_{y_i}(0)=x_i$, $q_{y_i}(1)=x'$ and $f\circ q_{y_i}=p_{y_i}$. If
$\bar q_{y_i}(t)=q_{y_i}(1-t)$ and $\bar p_{y_i}(t)=p_{y_i}(1-t)$,
we have $\bar q_{y_i}(0)=x'$, $f\circ \bar q_{y_i}=\bar p_{y_i}$.
Because $f(x')=\bar p_{y_i}(0)=y_i$, then $y_1=y_2$. Since the
lifting is unique (see for example Chapter 2 of \cite{spanier}), we
obtain that $\bar q_{y_1}=\bar q_{y_2}$; in particular $x_1=x_2$.
\end{proof}

Let $f:X\rightarrow Y$ be a covering projection between
path-connected metric spaces, and consider the associated morphism
between their fundamental groups $f_*:\pi_1(X)\rightarrow
\pi_1(Y)$. It is well known that $f$ is a homeomorphism onto $Y$
if and only if $f_*[\pi_1(X)]=\pi_1(Y)$ (see e.g.  \cite[Chapter
2]{spanier}). Thus we obtain at once the following corollary. In
particular,  if $X$ and $Y$ are Banach spaces, a local
homeomorphism has the continuation property for lines if and only
if it is a  global homeomorphism. In this way we obtain the
desired generalization of the result of Plastock \cite[Theorem
1.2]{plastock}.

\begin{corollary}\label{corolario0}
Let $f:X\rightarrow Y$ be a local homeomorphism between
path-connected metric spaces, where $Y$ is  $\mathcal P$-connected
and locally $\mathcal P $-contractible for some family $\mathcal P
$ of paths. Suppose that either $Y$ is simply connected or
$\pi_1(X)=\pi_1(Y)$ is finite. Then $f$ is a global homeomorphism
if and only if $f$ has the continuation property for every path in
$\mathcal P$.
\end{corollary}

Recall that a continuous map $f:X\rightarrow Y$ between topological
spaces is said to be {\it proper} if $f^{-1}(K)$ is a compact set in
$X$ whenever $K$ is a compact set in $Y$. More generally, we say
that $f$ is {\it weakly proper} if, for every compact subset $K$ of
$Y$, each connected component of $f^{-1}(K)$ is compact in $X$.

\begin{corollary}\label{corolario1}
Let $X$ and $Y$ be  metric spaces, and suppose that $Y$ is
$\mathcal P $-connected and locally
 $\mathcal P $-contractible for some family $\mathcal P $
of paths. Then every weakly proper local homeomorphism
$f:X\rightarrow Y$ is a covering projection.
\end{corollary}

\begin{proof}
Let $p:[0,1]\rightarrow Y$ be a path in $\mathcal P $, consider
$0<b\leq 1$, and suppose that $q:[0,b)\rightarrow X$ satisfies
$f\circ q=p$ over $[0,b)$. Then we have that $  { \rm Im
\hspace{0.02in}} q$ is relatively compact in $X$, since it is
contained in a connected component of set $f^{-1}(  { \rm Im
\hspace{0.02in}} p)$. If we now choose a sequence $\{t_{n}\}$ in
$[0,b)$ convergent to $b$, there exists a subsequence $(t_{n_{k}})$
such that $\{q(t_{n_{k}})\}$ is convergent in $X$.
\end{proof}

\section{Mean value theorem on metric spaces}

Let $X$ and $Y$ be metric spaces. If $f:X\rightarrow Y$ is a
continuous map and $x\in X$ is not an isolated point of $X$, we
define the {\it lower and upper scalar derivatives} of $f$ at $x$ by
    $$D_x^-f=\liminf_{z\rightarrow
    x}\frac{d(f(z),f(x))}{d(z,x)},\hspace{0.2in}
    D_x^+f=\limsup_{z\rightarrow x}\frac{d(f(z),f(x))}{d(z,x)},$$
where $z$ is restricted to points of $X$ different from  $x$.

\

 We introduce the above quantities motivated by the work of F. John
\cite{john} in a Banach space context. The next lemma follows
immediately from the definition:

\begin{lemma}\label{lema2}
Let $V$ and $W$ be open sets in the metric spaces $X$ and $Y$
respectively, and suppose that $g=f|_V:V\rightarrow W$ is a
homeomorphism. If $x\in V$ is not an isolated point and $y=f(x)\in
W$, we have that $D_y^+(g^{-1})=(D_x^-g)^{-1}$ and $
D_y^-(g^{-1})=(D_x^+g)^{-1}.$
\end{lemma}
 When $X$ and $Y$ are Banach spaces and $f:X\rightarrow Y$ is
 differentiable at $x$,
F.John obtained in \cite{john}   that $D^+_xf=\|f'(x)\|$ and, if in
addition $f'(x)$ is invertible, $D_x^-f=\|f'(x)^{-1}\|^{-1}$. The
same statement holds for smooth mappings between connected and
complete Riemannian manifolds:

\begin{example}\label{ejemploR}
{\it Let $f:M\rightarrow N$ be a $C^1$ map between connected and
complete  Riemannian  manifolds. Then, for every  $x\in M$ we have
that $D_x^+f=\|df(x)\|$. If in addition $df(x)\in
\isom\left(T_xM;T_{f(x)}N\right)$, then
$D_x^-f=\|[df(x)]^{-1}\|^{-1}.$}

\

\begin{proof}
Our proof will  work both in the finite-dimensional and
infinite-dimensional cases. Let $x\in M$ and consider
$\varepsilon>0$. Since the map $\Vert df(\cdot)\Vert$ is continuous
on $M$, there exists $r>0$ such that $\vert \, \Vert df(x) \Vert -
\Vert df(y)\Vert \, \vert <\varepsilon$ for every $y\in B_{2r}(x)$.
Now for each $z\in B_r(x)$ with $z\neq x$ there exists a $C^1$ path
$\gamma_{z,\varepsilon}:[0,1]\rightarrow M$ such that
$\ell(\gamma_{z,\varepsilon})\leq d(x,z)+ \min\{r; \varepsilon;
\varepsilon d(x,z)\}.$ Thus for every $y\in\im
\gamma_{z,\varepsilon}$ we have that $d(x,y)\leq \ell
(\gamma_{z,\varepsilon})\leq \min\{2r; (1+\varepsilon)d(x,z)\}$, and
$\Vert df(y)\Vert\leq \Vert df(x)\Vert +\varepsilon$. Then
$$
d(f(x),f(z))\leq \ell(f\circ\gamma_{z,\varepsilon})
=\int_0^1\|(f\circ \gamma_{z,\varepsilon})'(t)\|dt
$$
$$
\leq \sup\{(\Vert df(y)\Vert : y\in \text
{Im}\gamma_{z,\varepsilon}\} \cdot
\int_0^1\|(\gamma_{z,\varepsilon})'(t)\|dt $$
$$\leq (\Vert
df(x)\Vert +\varepsilon) (1+\varepsilon) d(x,z).
$$
In this way we obtain that $D_x^+f\leq \|df(x)\|.$

\

For the reverse inequality, consider $v\in T_xM$ with $\Vert v \Vert
=1$. Let $\gamma(t)=\exp_x(tv) $ be the unique geodesic in $M$ such
that $\gamma(0)=x$ and $\gamma'(0)=v$, and define
$\sigma(t)=\exp_{f(x)}^{-1}\circ f\circ\exp_x(tv)$.  Note that the
path $\sigma:(-r,r)\rightarrow T_{f(x)}N$ is defined for some small
interval $(-r,r)$ and satisfies $\sigma(0)=f(x)$ and
$\sigma'(0)=df(x)(v)$. For $\vert t \vert$ small enough , we have by
\cite[Theorem VIII.6.4]{lang} that
$$\vert t \vert =\Vert tv
\Vert=d(\exp_x(tv),x)=d(\gamma (t), x)
$$
and also
$$\Vert \sigma (t)
\Vert=d(\exp_{f(x)}(\sigma(t)), f(x))=d(f(\gamma(t)), f(x)).
$$
Therefore
$$
\Vert df(x)(v)\Vert =\Vert \sigma'(0)\Vert = \lim_{t\to
0}\frac{\Vert\sigma(t)\Vert}{\vert t \vert}= \lim_{t\to
0}\frac{d(f(\gamma(t)), f(x))}{d(\gamma(t), x)}\leq D_x^+f.$$

\

Suppose now that  $df(x)\in \isom\left(T_xM;T_{f(x)}N\right)$. Then
there exist $V$ and $W$,  open neighborhoods of $x$ and $f(x)$
respectively, such that $f|_V:V\rightarrow W$ is a diffeomorphism.
Using  Lemma \ref{lema2} we obtain that
$(D_x^-f)^{-1}=D_{f(x)}^+(f|_V)^{-1}=\Vert d(f|_V)^{-1}
(f(x))\Vert=\Vert [df(x)]^{-1}\Vert$.
\end{proof}
\end{example}

The first part of the proof in the above example also works in the
case of Finsler manifolds. Thus we obtain the following:

\begin{example}\label{ejemploF}
{\it Let $f:M\rightarrow N$ be a $C^1$ map between connected and
complete $C^1$ Finsler manifolds. Then, for every  $x\in M$ we have
that $D_x^+f \leq \|df(x)\|$. If in addition $df(x)\in
\isom\left(T_xM;T_{f(x)}N\right)$, then
$D_x^-f\geq\|[df(x)]^{-1}\|^{-1}.$}
\end{example}

\begin{remark}
{\it Comparison of $D_x^-f$ with the Ioffe-Katriel surjection
constant.} The surjection constant was introduced by Ioffe
\cite{ioffe} for non differentiable maps between Banach spaces to
establish global inversion theorems. Katriel \cite{katriel} also
works with the surjection constant in order to give global
homeomorphism theorems in certain metric spaces. It is defined as
follows. If $f:X\rightarrow Y$ a continuous mapping between metric
spaces we set, for $x\in X$ and $t>0$:
$${\rm Sur }(f,x)(t)=\sup\{r\geq 0: B_r(f(x))\subset f(B_t(x))\},$$
    $${\rm sur }(f,x)=\liminf_{t\rightarrow 0} ~ t^{-1}{\rm Sur }(f,x)(t).$$
Then ${\rm sur} (f,x)$ is called the {\it surjection constant of $f$
at $x$}. In general, ${\rm sur}(f,x)$  does not always coincide with
the lower scalar derivative $D_x^-f$. A simple example of this is
the inclusion map $i:\mathbb{R}\rightarrow\mathbb{R}^2$; it is easy
to calculate $D_0^-i=1$ and ${\rm sur}(i,0)=0$. Nevertheless, {\it
if $f:X\rightarrow Y$ is a local homeomorphism then, for all $x\in
X$
 $$D_x^-f={\rm sur}(f,x).$$}
\begin{proof}
Indeed, let $x\in X$ be fixed, and suppose that
$D_x^-f=\alpha<\gamma<\beta={\rm sur}(f,x)$. Since $f$ is a local
homeomorphism, there exists $\delta>0$ such that $f|_{B_\delta}(x)$
is a homeomorphism and ${\rm Sur}(f,x)(t)>t \gamma$ if $0<t<\delta$.
Therefore, $B_{t \gamma}(f(x))\subset f(B_t(x))$ if $0<t<\delta$.
Since $f$ is  injective on $B_\delta(x)$ we obtain for $0<t<\delta$
that if $d(x,z)=t$ then $f(z) \notin B_{t \gamma}(f(x))$ and
therefore $d(f(x),f(z))\geq t \gamma$.  In conclusion
$\alpha=D_x^-f\geq\gamma$, a contradiction. Thus ${\rm sur}(f,x)\leq
D_x^-f.$

Now suppose that $\beta={\rm sur}(f,x)<\gamma<\alpha=D_x^-f$. There
exist an open neighborhood $V^x$ and $r>0$ such that
$f|_{V^x}:V^x\rightarrow B_r(f(x))$ is a homeomorphism and
$d(f(x),f(z))>d(x,z)\cdot \gamma$, for all $z\in V^x$. Let $y\in Y$
such that $d(y,f(x))<\gamma\cdot t<r$. Then, there exists $z\in V^x$
with $y=f(z)$ and  $d(x,z)<t$. Therefore, ${\rm Sur}(f,x)(t)\geq
\gamma \cdot t$ if $0<t<\frac{r}{\gamma}$. Then $\beta={\rm
sur}(f,x)\geq\gamma$, a contradiction. Thus $D_x^-f\leq {\rm
sur}(f,x).$
\end{proof}
\end{remark}

In order to obtain  mean value inequalities  in terms of the lower
and upper scalar derivatives, we will prove two very simple, but
useful lemmas. If $q$ is a path defined on $[\alpha,\beta]$, we will
denote by $q_{t,s}$ the restriction of $q$ over
$[t,s]\subseteq[\alpha,\beta]$.

\begin{lemma}\label{lema3}
Let $f:X\rightarrow Y$ be a continuous map between metric spaces.
Suppose that $q:[\alpha,\beta]\rightarrow X$ is a rectifiable path
such that $\ell (q_{t,s})>0$ for all $t,s\in [\alpha,\beta]$ with
$t<s$, and denote $p=f\circ q$. Then,  for every
$t\in(\alpha,\beta)$ we have:
    $$\frac{d(p(\alpha),p(\beta))}{\ell(q_{\alpha,\beta})}\leq\max\Bigg{\{}
      \frac{d(p(\alpha),p(t))}{\ell(q_{\alpha,t})},
      \frac{d(p(t),p(\beta))}{\ell(q_{t,\beta})}
      \Bigg{\}}.$$
\end{lemma}

\begin{proof}
Let $t\in(\alpha,\beta)$ be fixed and suppose first that
$$
\ell(q_{t,\beta})d(p(\alpha),p(t))
\leq \ell(q_{\alpha,t})d(p(t),p(\beta)).
$$
Then
$$\ell(p_{t,\beta}) d(p(\alpha),p(\beta))
\leq \ell(p_{t,\beta}) [d(p(\alpha),p(t))+d(p(t),p(\beta))]
$$
$$
\leq \ell[(q_{\alpha,t})+ell(p_{t,\beta})] d(p(t),p(\beta))=
\ell(q_{\alpha,\beta}) d(p(t),p(\beta)).
$$
If on the other hand
$$\ell(q_{t,\beta})d(p(\alpha),p(t))
\geq \ell(q_{\alpha,t})d(p(t),p(\beta)),$$
we obtain in the same way
that
$$\ell(q_{\alpha, t}) d(p(\alpha),p(\beta)) \leq
\ell(q_{\alpha,\beta})d(p(\alpha),p(t)).$$
\end{proof}

With the same proof of Lemma \ref{lema3} we have the following
result:

\begin{lemma}\label{lema4}
Let $f:X\rightarrow Y$ be a continuous map between metric spaces.
Let $q:[\alpha,\beta]\rightarrow X$ be a  path such that
$q(\alpha)\neq q(\beta)$, and suppose that $p=f\circ q$ is
rectifiable. Then, for every $t\in(\alpha,\beta)$ we have:
    $$\frac{\ell(p_{\alpha,\beta})}{d(q(\alpha),q(\beta))}\geq\min\Bigg{\{}
      \frac{\ell(p_{\alpha,t})}{d(q(\alpha),q(t))},
      \frac{\ell(p_{t,\beta})}{d(q(t),q(\beta))}
      \Bigg{\}}.$$
\end{lemma}

\begin{remark}
Note that, under the hypothesis of Lemma \ref{lema4},
$d(q(\alpha),q(t))$ and $d(q(t),q(\beta))$ cannot be be $0$
simultaneously. If for example $t\in(\alpha,\beta)$ is such that
$q(t)=q(\alpha)$, we understand that
$\frac{\ell(p_{\alpha,t})}{d(q(\alpha),q(t))}:=\infty$. In this case
$q(t)\neq q(\beta)$ and  we have
$$
\frac{\ell(p_{\alpha,\beta})}{d(q(\alpha),q(\beta))}\geq
\frac{\ell(p_{t,\beta})}{d(q(t),q(\beta))}.
$$

\end{remark}

\begin{theorem}\label{teorema2}
Let $f:X\rightarrow Y$ be a continuous map between metric spaces,
let $q:[a,b]\rightarrow X$ be a rectifiable path, and denote
$p=f\circ q$. Then:
\begin{enumerate}
\item There exists $\tau\in [a,b]$ such that $d(p(a),p(b))\leq
D_{q(\tau)}^+f \cdot \ell(q).$
\item  We have $\ell(p) \leq \sup_{x\in \im q}D_x^+f \cdot \ell(q).$
\end{enumerate}
\end{theorem}

\begin{proof}
In order to prove part $(1)$, first note that using a suitable
reparametrization of $q$ (for example, the reparametrization by arc
length, see \cite[I.1.20]{bridson}) we may assume  with no lost of
generality that $\ell (q_{t,s})>0$ for all $t,s\in [a,b]$ with
$t<s$. Then we denote
$$\Delta_{t,s}f:=\frac{d(p(t),p(s))}{\ell(q_{t,s})}.$$
Let $\alpha_0=a$, $\beta_0=b$ and consider the midpoint $m_0\in
[\alpha_0,\beta_0]$. By  Lemma \ref{lema3} we have
    $$\Delta_{\alpha_0,\beta_0}f\leq
    \max\{\Delta_{\alpha_0,m_0}f,\Delta_{m_0,\beta_0}f\}:=\Delta_{\alpha_1,\beta_1}f.$$
Now take the midpoint $m_1\in[\alpha_1,\beta_1]$, and again by Lemma
\ref{lema3}, we get
    $$\Delta_{\alpha_1,\beta_1}f\leq
    \max\{\Delta_{\alpha_1,m_1}f,\Delta_{m_1,\beta_1}f\}:=\Delta_{\alpha_2,\beta_2}f.$$
By proceeding in this way, we construct a nested sequence of
intervals $[\alpha_n,\beta_n]$ with $\alpha_n<\beta_n$ and
$(\beta_n-\alpha_n)\rightarrow 0$, satisfying:
    $$\Delta_{a,b}f\leq \Delta_{\alpha_1,\beta_1}f
    \leq \cdots
    \leq\Delta_{\alpha_n,\beta_n}f\leq \cdots$$
Then there exists a point $\tau\in[a,b]$ such that
$\bigcap_n[\alpha_n,\beta_n]=\{\tau\}$. In the case that
$\alpha_n<\tau<\beta_n$ for every $n\in \mathbb N$,  we have
$$
\limsup_{n \to \infty}\Delta_{\alpha_n,\tau}f \leq \limsup_{n\to
\infty}\frac{d(f(q(\alpha_n)),f(q(\tau)))}{d(q(\alpha_n),q(\tau))}
\leq D_{q(\tau)}^+f
$$
and
$$
\limsup_{n \to \infty}\Delta_{\tau, \beta_n}f \leq \limsup_{n\to
\infty}\frac{d(f(q(\tau)),f(q(\beta)))}{d(q(\tau),q(\beta_n))} \leq
D_{q(\tau)}^+f
$$
so we obtain that
$$
\Delta_{a,b}f\leq \limsup_{n\to
\infty}\Delta_{\alpha_n,\beta_n}f\leq \limsup_{n\to\infty}
\max\{\Delta_{\alpha_n,\tau}f,\Delta_{\tau,\beta_n}f\} \leq
D_{q(\tau)}^+f.
$$
In the same way, we also obtain that $\Delta_{a,b}f\leq
D_{q(\tau)}^+f$ if either $\tau=\alpha_n$ for some $n$ or
$\tau=\beta_n$ for some $n$.

\

To prove part $(2)$, we may assume  that $K:=\sup_{x\in \im q}D_x^+f
<\infty$ (since otherwise the result holds trivially). Consider a
partition $t_0=a\leq t_1\leq ...\leq t_n=b$ of the interval $[a,b]$.
For each $i=1, \cdots,n$, applying part $(1)$ to the interval
$[t_{i-1},t_i]$ we have that
    $$d(p(t_{i-1}),p(t_i))\leq K \cdot
    \ell\left(q|_{[t_{i-1},t_i]}\right). $$
Therefore,
    $$\sum_{i=1}^{n}d(p(t_{i-1}),q(t_i))\leq
    K \sum_{i=1}^{n}
    \ell \left(q|_{[t_{i-1},t_i]}\right)=K\cdot \ell (q).$$
Taking  the supremum over all partitions of $[a,b]$ we conclude the
proof.
\end{proof}

Next we give an analogous result for the lower scalar derivative:

\begin{theorem}\label{teorema3}
Let $f:X\rightarrow Y$ be a continuous map between metric spaces,
let $q:[a,b]\rightarrow X$ be a path,  and suppose that $p=f\circ q$
is rectifiable. Then:
\begin{enumerate}
\item If $q(a)\neq q(b)$, there exists $\tau\in [a,b]$
 such that $\ell(p)\geq D_{q(\tau)}^-f \cdot d(q(a),q(b)).$
\item If $0< \inf_{x\in \im q}D_x^-f <\infty$, we have that
$\ell(p) \geq \inf_{x\in \im q}D_x^-f \cdot \ell(q).$
\end{enumerate}
\end{theorem}
\begin{proof}
Using Lemma \ref{lema4} we can construct a nested sequence of
intervals $[\alpha_n,\beta_n]\subset [a,b]$ with $\alpha_n<\beta_n$,
$g(\alpha_n)\neq g(\beta_n)$ and $(\beta_n-\alpha_n)\rightarrow 0$,
satisfying:
$$
\frac{\ell(p_{a,b})}{d(q(a),q(b))} \geq
\frac{\ell(p_{\alpha_1,\beta_1})}{d(q(\alpha_1),q(\beta_1))} \geq
\cdots \geq
\frac{\ell(p_{\alpha_1,\beta_1})}{d(q(\alpha_n),q(\beta_n))}\geq
\cdots
$$
Then we consider $\bigcap_n[\alpha_n,\beta_n]=\{\tau\}$, and  we
proceed as in the proof of Theorem \ref{teorema2}.
\end{proof}

\section{Covering projections via the bounded path-lifting property}

Let $f:X\rightarrow Y$ be a continuous map between path-connected
metric spaces, and let  $p:[0,1]\rightarrow Y$ be a path in $Y$.
We will say that $f$ has the {\it bounded path-lifting property
for $p$} if, for every $b\in (0,1]$ and every $q:[0,b)\rightarrow
X$ such that $f\circ q=p$ over $[0,b)$, there exists $\alpha>0$
such that:
$$\inf \{D_x^-f: x\in\im q\}\geq \alpha.$$

\begin{theorem}\label{teorema4}
Let $f:X\rightarrow Y$ a continuous map between metric spaces, and
suppose that $X$ is complete. If $f$ has the bounded path-lifting
property for rectifiable paths, then $f$ has the continuation
property for rectifiable paths. As a consequence, if   $Y$ is
path-connected and  locally  $\mathcal R$-contractible, and $f$ is
a local homeomorphism, then $f$ is a covering projection.
\end{theorem}

\begin{proof}
Let $p:[0,1]\rightarrow Y$ be a rectifiable path, and consider $b\in
(0,1]$ and $q:[0,b)\rightarrow X$ such that $f\circ q=p$ over
$[0,b)$. If $f$ has the bounded path-lifting property for $p$, there
exists $\alpha>0$ such that $\inf \{D_x^-f: x\in\im q\}\geq \alpha$.
By using Theorem \ref{teorema3} we obtain that, for every $s,t\in
[0,b)$ with $s<t$:
$$ \ell \left(p|_{[s,t]}\right)\geq  \alpha \cdot d(q(s),q(t)) .$$
Now let $\{t_{n}\}$ be an increasing sequence in $[0,b)$ convergent
to $b$. Then for all $m>n$, we have:
$$d(q(t_n),q(t_m))\leq \frac{1}{\alpha} ~
\ell\left(p|_{[t_n,t_m]}\right) \leq \frac{1}{\alpha} ~
\ell\left(p|_{[t_n,b]}\right).$$ Since $p$ is rectifiable, the map
$t\mapsto \ell\left(p|_{[t,b]}\right)$ is continuous (see, for
example, \cite[1.20(5)]{bridson}). This implies that
 $\{q(t_n)\}$ is a Cauchy
sequence, and therefore convergent since $X$ is complete. So, $f$
has the continuation property for $p$. Now the last part of the
result follows directly from Theorem \ref{teorema1}.
\end{proof}

Our next Corollary extends a classical result due to Ambrose for
smooth mappings between Riemannian manifolds (see \cite[Theorem
A]{ambrose} and see also \cite[Theorem VIII.6.9]{lang}).

\begin{corollary}\label{corolario2}
Let $f:X\rightarrow Y$ be a  local homeomorphism between
 metric spaces, where $X$ is complete and $Y$ is path-connected and
locally  $\mathcal R$-contractible. If there exists $\alpha>0$
such that $D_x^-f\geq \alpha$ for all $x\in X,$ then $f$ is a
covering projection.
\end{corollary}

Another easy consequence of Theorem \ref{teorema4} is the following.

\begin{corollary}\label{corolario3}
Let $f:X\rightarrow Y$ be a  local homeomorphism between metric
spaces, where $X$ is complete and $Y$ is path-connected and
locally  $\mathcal R$-contractible. Suppose  that:

\begin{enumerate}
\item For every bounded subset
$B$ of  $X$, we have  $\inf_{x\in B}(D_x^-f)>0$.

\item For some $y_0\in Y$ and  $x_0\in
X$, we have  $d(f(x),y_0)\rightarrow\infty$ as
$d(x,x_0)\rightarrow\infty$.
\end{enumerate}
Then, $f$ is a covering projection.
\end{corollary}

\begin{proof}
Let $p:[0,1]\rightarrow Y$ be a rectifiable path, let $b\in(0,1]$
and let $q:[0,b)\rightarrow X$ be such that $f\circ q=p$ over
$[0,b)$. Consider $R:=\max \{d(p(t),y_0):0\leq t \leq1\}$. There
exists $r>0$ such that $d(f(x),y_0)>R$ whenever $d(x,x_0)>r$. Then
    $$\inf\{D_x^-f:x\in\im q\}\geq
    \inf\{D_x^-f:d(f(x),y_0)\leq R\}$$
    $$
    \geq \inf\{D_x^-f:d(x,x_0)\leq r\}>0.\hspace{0.1in}^\Box$$
\end{proof}

The conditions of Corollaries \ref{corolario2} and \ref{corolario3}
appear frequently in global inversion theorems. For example,
Corollary \ref{corolario2} extends the classical well known result
in \cite[Theorem 1.22]{schwartz} concerning  $C^1$ mappings between
Banach spaces, and also extends  the analogous result of F. John
(see \cite{john}, Corollary in pg. 87) in the context of nonsmooth
mappings between Banach spaces. On the other hand, Corollary
\ref{corolario3} is an extension  of \cite[Corollary 3.3]{zampieri}
by Zampieri. Similar global inversion results for metric spaces
(with more complicated topological hypothesis over $X$ and $Y$) have
been obtained  by Katriel in \cite[Theorem 6.1]{katriel} and
\cite[Theorem 6.2]{katriel}.

\

Next we are going to see that the bounded path-lifting property
can be equivalently defined in terms of a weight. Here, be a {\it
weight} we mean  a nondecreasing map (not necessarily continuous)
$\omega:[0,\infty)\rightarrow (0,\infty)$ such that:
$$\int_0^\infty \frac{dt}{\omega(t)}=\infty.$$
Now let $f:X\rightarrow Y$ be a continuous map between
path-connected metric spaces. If $p:[0,1]\rightarrow Y$ is a path,
we will say that $f$ has the {\it bounded path-lifting property
for $p$ with respect to the weight $\omega$} if,  for every $b\in
(0,1]$ and $q:[0,b)\rightarrow X$ such that $f\circ q=p$ over
$[0,b)$, there exist $x_0\in X$ and $\alpha>0$ such that:
$$\inf \{D_x^-f\cdot  \omega (d(x,x_0)): x\in\im q\}\geq \alpha.$$

\begin{lemma}\label{proposicion1}
Let  $f:X\rightarrow Y$ be a continuous map between metric spaces,
and let $p:[0,1]\rightarrow Y$ be a rectifiable path in $Y$.  Then,
$f$ has the bounded path lifting property for $p$ if and only if $f$
has the bounded path lifting property for $p$ with respect to some
weight $\omega$.
\end{lemma}

\begin{proof}
The sufficiency follows trivially by choosing $\omega(t)\equiv 1$.
We are going to prove the necessity. Let $b\in (0,1]$ and
$q:[0,b)\rightarrow X$ such that $f\circ q=p$ over $[0,b)$. There
exist a weight $\omega$, some $x_0\in X$ and some $\alpha>0$ such
that $\inf \{D_x^-f\cdot  \omega (d(x,x_0)): x\in\im q\}\geq
\alpha$. With no loss of generality, we can suppose that
$x_0=q(0)$, since otherwise we could consider the alternative
weight $\bar \omega (t):=\omega(d(x_0,q(0))+t)$. Now  define the
map $\xi:[0,b)\rightarrow\mathbb{R}$ by
    $$\xi(t)=\max_{\tau\in[0,t]} d(q(\tau),x_0).$$
It is clear that $\xi$ is continuous and non-decreasing. Before
going further, we are going to show that:
    \begin{equation}\label{ecuacion}
\frac{\xi(t)-\xi(t')}{\omega(\xi(t))}\leq \frac{1}{\alpha} \cdot
\ell\left(p|_{[t',t]}\right),\hspace{0.3in} \forall\hspace{0.05in}
t' \leq t.
    \end{equation}
\noindent Indeed, let $t' <t$ in $[0,b)$, and consider the
interval $[t', t]$. Taking into account that
$\omega(d(q(\tau),x_0)\leq \omega(\xi(t))$ for every $\tau\in [t',
t]$, and using Theorem \ref{teorema3},  we obtain that
$$
\ell\left(p|_{[t',t]}\right) \geq
\frac{\alpha}{\omega(\xi(t))}\cdot d(q(t'),q(t)).
$$
Therefore
$$
d(q(t),x_0) \leq d(q(t'),x_0)+\frac{\omega(\xi(t))}{\alpha}
\cdot\ell\left(p|_{[t',t]}\right)\leq \xi(t')+
\frac{\omega(\xi(t))}{\alpha} \cdot\ell \left(p|_{[t',t]}\right).
$$
In order to establish \eqref{ecuacion}, note that the inequality
is clear if $\xi(t')=\xi(t)$. On the other hand, if
$\xi(t')<\xi(t)$, there exists  $t^*\in (t',t]$ such that
$\xi(t)=d(q(t^*), x_0)$. In this case, by applying the above
argument to $[t',t^*]$, we obtain that
$$\xi(t)=d(q(t^*), x_0)\leq \xi(t')+
\frac {\omega(\xi(t^*))}{\alpha}
\cdot\ell\left(p|_{[t',t^*]}\right)$$
$$\leq \xi(t')+\frac {\omega(\xi(t))}
{\alpha}\cdot\ell\left(p|_{[t',t]}\right).$$

Now let $0<\delta<b$ be fixed. Given a partition
$0=t_0<t_1<...<t_n=\xi(\delta)$ of $[0,\xi(\delta)]$, since $\xi$
is  continuous and non-decreasing we can find
 $0=t_0<t_1<\ldots<t_n=\delta$ such that $s_i=\xi(t_i)$, for
$i=0,\ldots,n$. Then, by inequality \eqref{ecuacion}, we have
$$
\sum_{i=1}^{n}\frac{s_i - s_{i-1}}{\omega(s_i)}
\leq\frac{1}{\alpha}\cdot\sum_{i=1}^{n}
\ell\left(p|_{[t_{i},t_{i-1}]}\right)= \frac{1}{\alpha}\cdot
\ell\left(p|_{[0,\delta]}\right)
\leq\frac{1}{\alpha}\cdot\ell(p).
$$
Therefore, for every $\delta \in [0,b)$ we obtain that
$$
\int_{0}^{\xi(\delta)}\frac{dt}{\omega(t)}\leq
\frac{1}{\alpha}\cdot\ell(p)<\infty.
$$
Since $\omega$ is a weight, we conclude that there exists some
$r>0$ such that $\xi(\delta)\leq r$ for every $\delta\in[0,b)$. As
a consequence, for every $x\in\im q$ we have that
$\omega(d(x,x_0))\leq\omega(r)$ and since $D_x^-f\cdot
\omega(d(x,x_0))\geq \alpha$ we finally have that
$$\inf \{D_x^-f:x\in\im q\}\geq
\frac{\alpha}{\omega(r)}>0.\hspace{0.1in}^\Box$$
\end{proof}

In order to define the Hadamard integral condition for a map
$f:X\rightarrow Y$ between  metric spaces, we need to restrict
ourselves to the case of mappings satisfying $0<D_x^-f<\infty$, for
every $x\in X$. In this case we say that $f$ is a {\it regular map}.
Now if $f:X\rightarrow Y$ is a regular between  metric spaces, we
will say that {\it $f$ satisfies the Hadamard integral condition}
if, for some $x_0\in X$,
$$\int _0^\infty \inf_{x\in \overline{B_t(x_0)}}D_x^-f ~
dt=\infty.$$

\begin{lemma}\label{observacion1}
Let $f:X\rightarrow Y$ be a regular map between metric spaces.
Then $f$ satisfies the Hadamard integral condition if and only if
there exist $x_0\in X$ and a weight $\omega$ such that
$D_x^-f\cdot\omega(d(x,x_0))\geq 1$, for every $x\in X$.
\end{lemma}
\begin{proof} Suppose that $f$ satisfies the Hadamard integral
condition for some $x_0\in X$. Then,  for every $t\geq 0$:
$${\inf\{ D_x^-f: x\in \overline{B_t(x_0)}\}}>0.$$
If we define $\omega(t)=[\inf\{ D_x^-f: x\in
\overline{B_t(x_0)}\}]^{-1}$, it is clear that  $\omega$ is a
weight and $D_x^-f\cdot\omega(d(x,x_0))\geq 1$, for all $x\in X$.

Conversely, suppose that there exist $x_0\in X$ and a weight
$\omega$, such that $D_x^-f\cdot\omega(d(x,x_0))\geq 1$, for all
$x\in X$. For every $t\geq 0$, and every $0\leq r\leq t$, we have
    $$\frac{1}{\omega(t)}\leq \frac{1}{\omega(r)}
    \leq\inf\{D_x^-f: d(x,x_0)=r\}.$$
Then $\omega(t)^{-1}\leq\inf\{D_x^-f:x\in \overline{B_t(x_0)}\}$
for all $t\geq 0$, and therefore $f$ satisfies the Hadamard
integral condition.
\end{proof}

Using Theorem \ref{teorema4} and Lemmas \ref{proposicion1} and
\ref{observacion1}, we deduce at once the following result. This
gives the desired extension of Hadamard Theorem to our context,
and extends also the analogous results of F. John (see
\cite{john}, Corollary in page 91) and Ioffe (see \cite[Theorem
2]{ioffe}).

\begin{theorem}\label{teorema5}
Let $f:X\rightarrow Y$ be a regular local homeomorphism between
metric spaces. Suppose that $X$ is complete and $Y$ is
path-connected locally $\mathcal R$-contractible. If $f$ satisfies
the Hadamard integral condition, then $f$ is a covering
projection.
\end{theorem}

As a direct application, using Example \ref{ejemploF}, we obtain
the following result (compare with \cite[Corollary 3.4]{GJ}).

\begin{corollary}
Let $f:M\rightarrow N$ a $C^1$ map between connected $C^1$  Finsler
manifolds, $M$ complete. Suppose that $df(x)$ is invertible, for
every $x\in M$. If
$$\int _0^\infty \inf_{x\in
\overline{B_t(x_0)}}\|[df(x)]^{-1}\|^{-1} ~ dt=\infty,$$ for some
$x_0\in M$, then $f$ is a covering projection.
\end{corollary}

To finish this section, we note that  Hadamard condition is not
necessary for a local homeomorphism to be a covering projection,
even in very simple cases. For instance, the map $f(x,y)=(x+y^3,y)$
is a global homeomorphism from $\mathbb R^2$ to $\mathbb R^2$, but
it does not satisfy the Hadamard integral condition (see
\cite[Example 1.7]{xavier}).

\section{Locally quasi-isometric maps}

In this section, we obtain a more complete result that relates the
lifting concepts that we use in this paper, but we need an extra
assumption on regularity. We shall say that a map $f:X\rightarrow Y$
between metric spaces is a {\it quasi-isometry}  if $f$ is a
homeomorphism and there exist $0<\alpha\leq \beta<\infty$ such that:
    $$\alpha\leq\inf_{x\in X}D_x^-f\leq\sup_{x\in X}D_v^+f\leq \beta.$$
In the same way, we will say that  $f:X\rightarrow Y$  is a {\it
local quasi-isometry}  if for every $x\in X$ there exist open
neighborhoods $V$ of $x$ and  $W$ of $f(x)$ such that
$f|_V:V\rightarrow W$ is a quasi-isometry.

\

It is clear that if a map between metric spaces is locally a
bi-Lipschitz homeomorphism, then it is a local quasi-isometry.
More generally, using Example \ref{ejemploF} we obtain that every
local $C^1$-diffeomorphism between Banach spaces or between
Finsler manifolds is also a local quasi-isometry.

\

Recall that, for a metric space $X$,  the {\it length} of a path
$p:[a,b)\rightarrow X$ defined on a semi-open interval is defined
by:
    $$\ell(p):=\lim_{t\rightarrow b^-}\ell\left(p|_{[a,t]}\right).$$
In this case we also say that the path $p$ is  {\it rectifiable}
when $\ell(p)<\infty$. We will need the following simple Lemma.

\begin{lemma}\label{lema5}
Let  $q:[a,b)\rightarrow X$ be a rectifiable path on a metric space.
Then for every  sequence $\{t_n\}\subset [a,b)$ converging to $ b$,
the sequence $\{q(t_n)\}$ is a Cauchy sequence in $X$.
\end{lemma}

\begin{proof}
Consider a sequence $\{t_n\}\subset [a,b)$ converging to $b$. For
each $\varepsilon>0$ there exists $n_0\in \mathbb N$ such that
$\ell(p) - \ell\left(q|_{[a,t_n]}\right)< \varepsilon$ for every
$n\geq n_0$. If $n,m\geq n_0$ and $t_n\leq t_m$, then:
$$d(q(t_n), q(t_m))\leq \ell\left(q|_{[t_n,t_m]}\right)
=\ell\left(q|_{[a,t_m]}\right)-\ell\left(q|_{[a,t_n]}\right)
<\varepsilon.$$\hfill\end{proof}

Now we can easily derive our main result in this Section.  This
extends \cite[Theorem 3.5]{GJ}. We note that Condition 3 was the
key of the original argument  used by Hadamard in \cite{hadamard}.
On the other hand, Condition (8) was introduced by Rabier
\cite[Theorem 5.3]{rabier} in the context of Finsler manifolds.

\begin{theorem}\label{teorema6}
Let  $f:X\rightarrow Y$ be a local quasi-isometry  between
complete metric spaces, and suppose that $Y$ is path-connected and
locally $\mathcal R $-contractible. Then the following statements
are equivalent:
\begin{enumerate}
  \item $f$ is a covering projection.
  \item $f$ has the continuation property for rectifiable paths.
  \item For every path $q:[a,b)\rightarrow X$,
  we have $\ell(q)<\infty$ whenever $\ell (f\circ q)<\infty$.
  \item $f$ has the bounded-path lifting property for rectifiable
  paths.
  \item $f$ has the bounded-path lifting property for rectifiable
  paths with respect to some weight.
\end{enumerate}
\noindent If in addition, we assume that either $Y$ is simply
connected or $\pi_1(X)=\pi_1(Y)$ is finite, then the previous
conditions are also equivalent to the following:
\begin{enumerate}
\setcounter{enumi}{5}
\item $f$ is a homeomorphism.
\item $f$ is a proper map.
\item For every compact subset $K\subset Y$, there is a constant
$\alpha_K>0$ such that $D_x^{-}f \geq \alpha_K$, for every $x\in
f^{-1}(K)$.
\end{enumerate}
\end{theorem}

\begin{proof}
$(5) \Leftrightarrow (4)$ It is given in Lemma \ref{proposicion1}.

\noindent $(4)\Rightarrow (3)$  Let $q:[a,b)\rightarrow X$ be a
path in $X$, and suppose that $\ell(f\circ q)<\infty$. By Lemma
\ref{lema5} we have that, for every  sequence $\{t_n\}\in[a,b)$
converging to $b$, the sequence $\{f\circ q ~ (t_n)\}$ is
convergent in $Y$. Then there exists a continuous path
$p:[a,b]\rightarrow Y$ such that $f\circ q=p$ over $[a,b)$. In
addition  $p$ is rectifiable since
$$\ell(p)=\lim_{t \to b^-} \ell \left(p|_{[a,t]}\right)=
 \ell(f\circ q)<\infty.$$
Therefore there exists $\alpha>0$ such that $\inf\{D_x^-f: x\in \im
q\}\geq \alpha$. Thus by Theorem \ref{teorema3} we get $\ell
(q)<{\ell (p)}\cdot{\alpha}^{-1}<\infty.$

\noindent $(3)\Rightarrow (2)$ Let $p:[0,1]\rightarrow Y$ be a
rectifiable path, $b\in (0,1]$ and $q:[0,b)\rightarrow X$ such that
$f\circ q=p$ over $[0,b)$.  Since $\ell
\left(p|_{[0,b]}\right)<\infty$, by the hypothesis we obtain that
$\ell(q)<\infty$. Therefore, by Lemma \ref{lema5}, for every
sequence $\{t_n\} \in[0,b)$ converging to $b$ we have that
$\{q(t_n)\}$ is convergent in $X$.

\noindent $(2)\Rightarrow (1)$ It follows from Theorem
\ref{teorema1}.

\noindent $(1)\Rightarrow (4)$ Let $p:[0,1]\rightarrow Y$ be a
rectifiable path, $b\in (0,1]$ and $q:[0,b)\rightarrow X$ such that
$f\circ q=p$ over $[0,b)$. If $f$ is a covering projection then $f$
lifts paths and has  the unique-path-lifting property (see for
example \cite[Section 2.2]{spanier}). Therefore $q$ can be
continuously extended to $[0,1]$, and in particular the closure of
$\im q$ is compact in $X$. Now since  $f$ is a local quasi-isometry,
a simple compactness argument gives that  there exists $\alpha>0$
such that $\inf\{D_x^-f: x\in\im q\}\geq \alpha.$

\noindent Finally, in the case that either $Y$ is simply connected
or $\pi_1(X)=\pi_1(X)$ is finite, it is clear  that $(1)
\Leftrightarrow (6) \Rightarrow (7)\Rightarrow (8)\Rightarrow (4)$.
\end{proof}

\section{Applications  to global implicit mapping theorems}

In this section we obtain general global implicit mapping theorems
for metric spaces using Theorem \ref{teorema1}. In the Banach space
case, we derive a global implicit mapping theorem in terms of
partial derivatives, using a growth condition  analogous to the
classical Hadamard integral condition.

Let $X$, $Y$ and $W$ be  metric spaces, consider a continuous map
$f:X\times Y\rightarrow W$, and let $w\in W$ be fixed. We denote
$$Z_w:=\{(x,y)\in X\times Y:f(x,y)=w\}.$$
We consider  $Z_w$ as a metric space with a metric inherited from
$X\times Y$.

We will say that there exist {\it local implicit maps for $f$ at
$w$} if for every $(x_0,y_0)\in Z_w$ there exist neighborhoods $U$
of $x_0$ and $V$ of
 $y_0$ and a {\it unique} continuous map $g:U\rightarrow V$ such that:
    $$(x,y)
    \in U\times V\hspace{0.05in}{\rm and}\hspace{0.05in} f(x,y)=w
    \iff x\in U\hspace{0.05in}{\rm and}\hspace{0.05in}y=g(x).$$
We will say that there exists a {\it global implicit map for $f$
at $w$} if there exists a  {\it unique} map $g:X\rightarrow Y$
such that
  $$(x,y)
    \in X\times Y\hspace{0.05in}{\rm and}\hspace{0.05in}f(x,y)=w
    \iff x\in X\hspace{0.05in}{\rm and}\hspace{0.05in}y=g(x).$$
Consider the natural projection $\pi:Z_w\rightarrow X$ given by
$\pi(x,y)=x$. It is easy to see that  there exist  local implicit
maps for $f$ at $w$ (respectively, a global implicit map) if, and
only if, $\pi$ is a local  homeomorphism (respectively, a global
homeomorphism). Bearing  this in mind, the following results can be
considered as  global implicit mapping theorems.

\begin{corollary}\label{coro7}
Let $X$, $Y$ and $W$ be  metric spaces, where $X$ is  simply
connected, $\mathcal P$-connected and and locally  $\mathcal
P$-contractible for some family $\mathcal P$ of paths. Suppose
that $\pi:{Z_w}\rightarrow X$ is a local homeomorphism satisfying
the continuation property for every path in $\mathcal P$. Then
$\pi:{Z_w}\rightarrow X$ is a covering projection and, for every
connected component $C_w$ of $Z_w$, we have that
$\pi|_{C_w}:{C_w}\rightarrow X$ is a homeomorphism.
\end{corollary}

\begin{proof}
From Theorem \ref{teorema1} we obtain that $\pi:{Z_w}\rightarrow X$
is a covering projection. On the other hand, since
$\pi:{Z_w}\rightarrow X$ is a local homeomorphism and $X$ is
path-connected, we have that $Z_w$ is locally path-connected.
Therefore if $C_w$ is a connected component of $Z_w$ we have that
$C_w$ is in fact path-connected. Then  as in
 Corollary \ref{corolario0} we obtain that
 $\pi|_{C_w}:{C_w}\rightarrow X$ is a homeomorphism.
\end{proof}

Note that if $X$ is a Banach space, Corollary \ref{coro7} applies to
the family $\mathcal L$  of all lines in $X$. Therefore, Corollary
\ref{coro7} generalizes \cite[Theorem 4]{blot}. Next we give a
direct application of  Corollary \ref{corolario1}:

\begin{corollary}\label{coro8}
Let $X$, $Y$ and $W$ be  metric spaces, where $X$ is simply
connected, $\mathcal P$-connected and and locally  $\mathcal
P$-contractible for some family $\mathcal P$ of paths. Suppose
that $Z_w$ is connected and $\pi:{Z_w}\rightarrow X$ is a local
homeomorphism. Then the following conditions are equivalent:
\begin{enumerate}
\item $\pi:{Z_w}\rightarrow X$ satisfies the continuation
property for every path in $\mathcal P$.
\item $\pi:{Z_w}\rightarrow X$ is a proper map.
\item $\pi:{Z_w}\rightarrow X$ is a homeomorphism.
\end{enumerate}
\end{corollary}

Corollary \ref{coro8}
 should be compared with  \cite[Lemma 1]{ichiraku}.
 Note that, under the hypothesis of Corollary \ref{coro8},  $Z_w$ is
connected if, and only if, for some $x\in X$, $\pi^{-1}(x)$ contains
only one point. On the other hand, the map $\pi:{Z_w}\rightarrow X$
is proper if and only if, whenever  $\{(x_n,y_n)\}$ is a sequence
contained in $Z_w$ such that $(x_n)$ is convergent in $X$, then
$(y_n)$ has a convergent subsequence in $Y$.

\

A direct application  of  Corollary \ref{corolario2} is the
following:

\begin{corollary}
Let $X$, $Y$ and $W$ be  metric spaces, where $X$ is simply
connected, path-connected and and locally  $\mathcal
R$-contractible. Suppose that $\pi:{Z_w}\rightarrow X$ is a local
homeomorphism, and there exists $C>0$ such that every
$(x_0,y_0)\in Z_w$ has a neighborhood $U$ in $Z_w$ satisfying that
$d(y,y_0)\leq C d(x,x_0)$ for every $(x,y)\in U$. Then
$\pi:{Z_w}\rightarrow X$ is a covering projection.
\end{corollary}
\begin{proof}
For every $(x,y)\in U$ with $(x,y)\neq(x_0,y_0) $ we have
$$
\frac{d(x,x_0)}{d(x,x_0)+d(y,y_0)}\geq \frac{1}{1+C}
$$
and therefore $\displaystyle{D^-_{(x_0,y_0)}\pi \geq
\frac{1}{1+C}>0}$.
\end{proof}
\

Our  last result generalizes  \cite[Theorem 5]{ichiraku} by
Ichiraku, who considered  finite dimensional  spaces and constant
weight.

\begin{corollary}\label{ichiripollo}
Let $E$, $F$ and $W$ be  Banach spaces, let $f:E\times
F\rightarrow W$ be a $C^1$ map and consider the set:
    $$Z_0:=\{(x,y)\in E\times F:f(x,y)=0\}.$$
Suppose that $\partial_y f(x,y)\in\isom(F;W)$ for all $(x,y)\in Z_0$
and that for some continuous weight $\omega:[0,\infty)\rightarrow
(0,\infty)$:
    $$\|\partial_yf(x,y)^{-1}\|\cdot\|\partial_xf(x,y)\|\leq\omega(\|y\|),
    \hspace{0.3in}\forall (x,y)\in Z_0.$$
If  $Z_0$ is connected, then there exists a global implicit map
for $f$.
\end{corollary}

\begin{proof}
Since $\partial_y f(x,y)\in\isom(F;F)$ for all $(x,y)\in Z_0$, there
exists a local implicit map for $f$, that is, $\pi:{Z_0}\rightarrow
E$ is a local homeomorphism. It is enough to show that
$\pi{Z_w}\rightarrow E$ has the continuation property for every
line. Let $p$ be a line in $E$, $b\in (0,1]$ and $y:[0,b)\rightarrow
F$ such that $(p(t),y(t))\in Z_0$ en $[0,b)$. Since $p$ is a line,
$p'(t)$ is constant, call it $v$. Therefore,
    $$y'(t)=-\partial_yf(p(t),y(t))^{-1}\partial_xf(p(t),y(t))v.$$
On the other hand,
    $$\|y(b)-y(0)\|\leq\ell \left(y|_{[0,b]}\right)=\int_0^b\| y'(t)\|dt\leq \|v\|\int_0^b\omega(\|y(t)\|)dt.$$
By \cite[Lemma 2.1]{radulescu}, we have
    $$\int_{\|y(0)\|}^{\|y(s)\|}\frac{dt}{\omega(t)}\leq\int_0^s\|v\|dt=\|v\|s\leq \|v\|,\hspace{0.3in}\forall s\in[0,b)$$
And, since $\omega$ is a weight, $\sup\{\|y(t)\|:t\in
[0,b)\}<\infty$. Because $\omega$ is nondecreasing, there exists
$\alpha>0$ such that $\omega(\|y(t)\|)\leq\alpha$, for every $t\in
[0,b)$. Then, for $s\in [0,b)$,
    $$\ell \left(y|_{[0,s]}\right)\leq \|v\|\int_0^s\omega(\|y(t)\|)dt\leq \|v\|\int_0^s\alpha dt\leq \|v\|\alpha s\leq \|v\|\alpha.$$
Therefore $\ell (y)<\infty$. By Lemma \ref{lema5} and the
completeness of $F$,  $\pi$  has the conti\-nuation property for
lines.
\end{proof}

\centerline{\sc Acknowledgements}

It is a great pleasure to thank Professor Gilles Godefroy for some
valuable conversations concerning the topic of this paper.



\end{document}